\documentclass[12pt]{amsart}
\usepackage{amsmath}

\newcommand{\E}{\overline{E}}
\newcommand{\C}{\mathbb C}

\newcommand{\Pn}{{\mathcal P}_n}

\newtheorem{thm}{Theorem}

\newtheorem{defi}{Definition}
\begin{document}
\title{A Proof of Sendov's Conjecture}
\author{Gerald Schmieder}
\date{}
\begin{abstract} The Theorem of Gau\ss{}-Lucas states that the
zeros of the derivative of a complex polynomial $p$ are contained
in the closed convex hull of its zeros.
For a polynomial $p$ having all its zeros in the closed unit disk Bl. Sendov
conjectured that the distance of an arbitrary zero to the closest derivative
zero is at most $1$. In this article we will give a proof.
\end{abstract}
\maketitle
\hspace{3mm}The zeros of the derivative of a complex polynomial
$p$ are functions of the zeros of $p$ itself. In general we do not
know explicit expressions for these functions. So approximate
localizations of the derivative zeros in terms of the given zeros
of $p$ are of interest.
A question of this type is the famous conjecture of Bl. Sendov
which goes back to 1959 and took place in Hayman's booklet on problems in
Complex Analysis (1967, \cite{hay},
by a misunderstanding, there it was named after Ilieff).
This conjecture states:\newline
\smallskip
{\em Let $p \in \C\, [z]$ be a polynomial of degree  $n>1$ having all
zeros
$z_1,\dots,z_n$ in the closed unit disk $\E$.
Does there exist for every $z_j$ some $\zeta$ with
$|z_j-\zeta|\le 1$ and $p'(\zeta)=0$ ?}\newline
\smallskip
For a history of the conjecture and a list of the numerous (about
100) publications on it, most of them in famous international
journals, see the recent article of Bl. Sendov \cite{sen}. In this
paper we give a proof of this question.
\newline
By $\Pn$ we denote the class of all monic polynomials of degree
$n>1$ having all its zeros in $\E$.\par For
the following we fix some polynomial $p\in\Pn$ with the
factorization $$p(z)=\prod _{j=1}^n (z-z_j).$$
\begin{defi}
Let $p\in\Pn$ and $w_0\in\C$ a zero of $p$. The disk $|z-w_0|\le
\rho$ is called critical with respect to $w_0$ if $p'$ has no zero
in the open disk but at least one on the boundary (the critical
circle). In this case $\rho=\rho(p,w_0)\ge 0$ is called the
critical radius for $w_0$ and the derivative zeros of $p$ on the
critical circle are called to be essential (with respect to
$w_0$). The polynomial $p\in\Pn$ is maximal with respect to the
point $w_0\in \E$ if among all polynomials $q\in\Pn$ with
$0=q(w_1)$ the critical radii fulfill $\rho(p,w_0)\ge
\rho(q,w_1)$.
\end{defi}
Of course the term $\rho(p,w_0)$ also makes sense for polynomials
$p$ with $p(w_0)=0$, which are not necessarily in $\Pn$. But in
this general case one will not find a maximal polynomial as this
is true in the compact class $\Pn$. \par Now Sendov's conjecture
can be formulated as $\rho(p,w_0)\le 1$ for all $p\in\Pn,\,
p(w_0)=0$. In order to prove the conjecture it would be enough to
check maximal polynomials in $\Pn$. But which $p\in\Pn$ are
maximal? Phelps and Rodriguez \cite{phro} guessed that these are
only the polynomials $p_n(z)=z^n-1$ and their rotations
$p_n(ze^{i\alpha})e^{-in\alpha}$. In the following we will confirm
this extension of Sendov's conjecture.
\par This will come out as a consequence of
\begin{thm}\label{thm1} Let $p\in \Pn$ have the zero  $z_1\in E$.
Then there is some $p^*\in\Pn$ which has a zero $w_0$ on the unit
circle and fulfills $\rho(p,z_1)\le \rho(p^*,w_0)$.\end{thm}

\section{The basic idea\label{basic}}
We start with some elementary formulas.
If $p\in\Pn$ is a polynomial with the zeros $z_1,\dots,z_n$ and
the derivative zero $\zeta$ with $p(\zeta)\neq 0$, then
\begin{eqnarray}\label{form1}
\frac{p'}{p}(\zeta)=0=\sum_{j=1}^n\frac{1}{\zeta-
z_j}.\end{eqnarray} We let the zeros $z_2,\dots,z_n$ be fixed and
vary $z_1$, i.e., we consider the polynomials
\begin{eqnarray}\label{Q}
Q(z,u)=(z-u)\,\prod _{j=2}^n(z-z_j)=(z-u)\,q(z).\end{eqnarray}
\smallskip
We assume for the moment that $\zeta$ is a zero of $p'$, but not a
zero of $p''$. The implicit function theorem (cf. \cite{nar})
shows the existence of a holomorphic function $\zeta(u)$ with
$\zeta(z_1)=\zeta$ and $\displaystyle \frac{\partial Q}{\partial
z}(\zeta(u),u)\equiv 0$, defined in a neighborhood of $z_1$. If we
move $u$ along a path $\gamma$ in $\C$ starting in $\gamma(0)=z_1$
then we have an unrestricted analytic continuation of
$\zeta(\gamma(t))$ if $\frac{\partial^2Q}{\partial
z^2}(\zeta(\gamma(t)),\gamma(t))\neq 0$ for all $t$. If the path
would meet these exceptional points, we would have at least a
continuation of $\zeta(\gamma(t))$ which is at least continuous in
such points. Note that the values of $\zeta(\gamma(t))$ with
respect to this continuation move on the Riemann surface $R$,
which is defined by the equation $Q'(z,u)=0$ (derivative with
respect to $z$). We will discuss this surface in section 2.
\smallskip
Note that $\ln |\zeta(u)-u|=\Re \log(\zeta(u)-u)$.
\par Let $p\in\Pn$ and $\zeta$ be a (not necessarily essential) derivative zero of $p$.
As above let $z_1,z_2,\dots,z_n\in\overline{E}$ be the zeros of
$p$ and $|z_1|<1$. If $\gamma:[0,1]\to \C$ is a path with
$\gamma(0)=z_1, \gamma(1)=u\in\C$, we see $$\frac{d}{dt}\ln
|\zeta(\gamma(t))-\gamma(t)|=\frac{d}{dt}\Re \log
(\zeta(\gamma(t))-\gamma(t))=\Re
\frac{\frac{d}{dt}(\zeta(\gamma(t))-\gamma(t))}
{(\zeta(\gamma(t))-\gamma(t))}.$$ Note that $\zeta(\gamma(t))$
depends on the path $\gamma$. Again we assume that
$\zeta(\gamma(0))=\zeta$. So we have $$\ln |\zeta(u)-u|-\ln
|\zeta-z_1|=\int_0^{\,1} \frac{d}{dt}\ln
|\zeta(\gamma(t))-\gamma(t)|\,dt$$$$=\Re\int_0^{\,1}
\frac{d}{dt}\log
(\zeta(\gamma(t))-\gamma(t))\,dt=\Re\int_0^{\,1}\gamma'(t)\frac{\zeta'(\gamma(t))-1}{\zeta(\gamma(t))-\gamma(t)}\,dt
=\Re\int_\gamma \frac{\zeta'(v)-1} {\zeta(v)-v}\,dv.$$ The right
hand side can be written as $$\Re \int_\gamma \frac{\zeta'(u)}
{\zeta(v)-v}-\frac{1}{\zeta(v)-v} \,dv.$$ From (\ref{form1}) we
obtain
$$0=\frac{Q'(\zeta(v),v)}{Q(\zeta(v),v)}=\frac{1}{\zeta(v)-v}+\frac{q'}{q}(\zeta(v)).$$
So it comes out
\begin{eqnarray}\nonumber
\ln |\zeta(u)-u|=\ln |\zeta-z_1|+\Re \int_\gamma \frac{\zeta'(u)}
{\zeta(v)-v}+ \frac{q'}{q}(\zeta(v))\,dv\,,
\end{eqnarray}
and therefore
\begin{eqnarray}\label{int}
|\zeta(u)-u|=|\zeta-z_1|\cdot\Big|\exp\left(\int_\gamma
\frac{\zeta'(u)} {\zeta(v)-v}+ \frac{q'}{q}(\zeta(v))\,dv
\right)\Big|.
\end{eqnarray}

\section{The Riemann surface $R$\label{Riemann}}
It is enough to consider polynomials $p\in Pn$ with the property
that $p$ has no multiple zeros and no multiple derivative zeros.
If we succeed to prove theorem \ref{thm1} under this restriction,
the general statement is clear by a continuity argument. By the
same argument we may assume that $q$ has no multiple zeros and no
multiple derivative zeros.\par The Riemann surface $R$ of the
derivative zeros of $Q$ is given by the equation
\begin{eqnarray}\label{gl} Q'(w)=q(w)+(w-u)q'(w)=0.\end{eqnarray}
This (actually compact) manifold $R$ consists of the points $w$
(which are the derivative zeros of $Q(.,u)$, and the equation
gives local uniformizations of $R$, if the derivative of
$u=\varphi(w):=w+\frac{q}{q'}(w)$ with respect to $w$ does not
vanish (note that these branch points are also described by
$\frac{\partial^2Q}{\partial z^2}(w,u)= 0$). So the points $w$
where $2q'(w)^2=q(w)q''(w)$ are branch points of the surface. this
branch points play in fact no special role on the Riemann surface,
their appearance depend on the special local coordinates, which
are given by the defining equation (example: the surface of the
square root is defined by $w^2=u$ with $0$ as a branch point; if
we add this point, it is conformally equivalent to the plane resp.
$\overline{\C}$). They can actually added as ''normal'' points to
the surface and have simply connected neighborhoods on which local
coordinates can be found.\par $R$, as a compact surface, may be
regarded as a $(n-1)$-sheeted covering of $\overline{\C}$, and
$\varphi$ gives a canonically projection $R\to \overline{\C}$.\par
We define
\begin{eqnarray}\label{deff}
f(u,\zeta(u)):=\exp\left(\int_\gamma \frac{\zeta'(u)}
{\zeta(v)-v}+ \frac{q'}{q}(\zeta(v))\,dv \right),\end{eqnarray}
where $\gamma_u:[0,1]\to \C$ with $\gamma_u,(0)=z_1,
\gamma_u(1)=u\in\C$ and $\zeta(\gamma_u(0))=\zeta$ (some fixed
derivative zero of $p$), $\zeta(\gamma_u(1))=\zeta(u)$. By
(\ref{int}) we have
\begin{eqnarray}\label{Qf}
|\zeta(u)-u|=|\zeta-z_1|\cdot |f(u,\zeta(u))|.\end{eqnarray} $f$
is, up to isolated singularities, a holomorphic function on $R$,
because it has this property in the local coordinate $u\in\C$ (the
case $u=\infty$ we discuss separately) . The holomorphy is not
obviously clear in the following cases.
\begin{enumerate}
\item [(i)] $Q(w_1,u_1)=0$ (this includes the case $u=\zeta(u)$), or
\item [(ii)] $2q'(w_2)^2=q(w_2)\cdot q''(w_2)$ (branch points)
\end{enumerate}
We discuss this two cases.\par\vspace{3mm} Case (i): The
assumption implies that $Q(.,u_1)$ has a multiple zero in the
point $w_1$. This is only possible if $u_1$ is one of the zeros
$z_1,\dots,z_{n-1}$ of $p$  and $u_1=w_1$. We have
$|\zeta(u_1)-u_1|=0$ if $\zeta(u_1)=w_1$. So this singularities of
$f$ is removable by (\ref{Qf}). Moreover we have
$\rho(Q(.,u_1))=0$ in this case.\par\vspace{3mm} Case (ii): If
$2q'(w_2)^2=q(w_2)q''(w_2)$, then $w_2\notin
\{z_1,z_2,\dots,z_{n-1}\}$, because $q$ has only simple zeros in
the points $z_2,\dots,z_n$. Especially $\frac{q'}{q}(w_2)$ is
finite. So (\ref{deff}) shows that $f$ is bounded in a
neighborhood of the branch point $w_2$ on $R$. Again we conclude
that $f$ has a removable singularity in this case.\medskip

We summarize: The function $f$ as defined in (\ref{deff}) is
holomorphic on the Riemann surface $R':=\{w\in
R\,:\,\varphi(w)\in\C\}$.
\medskip

From (\ref{Qf}) we obtain that $f(u,\zeta(u))$ equals
$\frac{\zeta(u)-u}{\zeta-z_1}$, up to a possible factor of modulus
one. For $u=z_1$ we see that this factor is one. By (\ref{gl}) we
receive the representation:
\begin{eqnarray}\label{fneu}
f(u,\zeta(u))= \frac{q'}{q}(\zeta) \cdot\frac{q}{q'}(\zeta(u)).
\end{eqnarray}
\par Finally we investigate the structure of $R$ close to $u=\infty$.
The point infinity is no branch point of $R$, because the function
$1/\varphi(1/w)$ has in $w=0$ the expansion
$w(\frac{n-1}{n}+a_1w+\dots)$. For $u\in \overline{E}$ all zeros
of $Q(.,u)$ are contained in $\overline{E}$. By the Gau\ss{}-Lucas
theorem we know that the zeros of the derivative
$Q'(z,u)=\frac{\partial Q}{\partial z}(z,u)$ lie in the convex
hull $C$ of the zeros. They are inner points of $C$ with the only
exception of multiple zeros of $Q$. None of these derivative zeros
in our case is of bigger order than $1$. So the same argument
gives that the zeros of the second order derivative
$Q''(z,u)=\frac{\partial^2Q}{\partial z^2}(z,u)$ are points the
open unit disk $E$. So the same is true for the branch points of
$R$. To be more precise, all branch points $w$ of $R$ fulfill
$|\varphi(w)|<1$.\par The subset $D_1$ of $R$ with
$\varphi(D_1)=\overline{E}$ therefore contains all branch
points.\par As a consequence, the complement $R\setminus D_1$
(including $\infty$) consists of $n-1$ simply connected domains
$G_1,\dots,G_{n-1}$. Let $\zeta(u)$ be the function which is
defined on $G_k$ with respect to a fixed start point $\zeta_0$
with $\varphi(\zeta_0)=z_1$. Then the mappings
$\Phi_k:=\varphi|G_k=\varphi|G_k:G_k\to \{u\in
\overline{\C}\,:\,|u|>1\}$ are conformal.\par The boundaries of
the domains $G_j$ are pairwise disjoint. Each $\partial G_j$ is
mapped homeomorphically by $\varphi$ on the unit circle.\par It
holds $P(z,u):=\frac{Q(z,u)}{u}=(\frac{z}{u}-1)q(z)$. The
derivative zeros of $P$ with respect to $z$ are the same as those
of $Q$. For $u\to\infty$ the polynomials $P(z,u)$ tend locally
uniformly to $q(z)$. So, in this case, $\zeta(u)$ tends to
$\infty$ on one $G_k$, let us say on $G_1$. For $k=2,\dots,n-1$ it
follows that each $\zeta_j(u)\in G_k$ tends to some derivative
zero $\xi_k$ of $q'$ if $u\to \infty$. \par If $\zeta(u)\in G_1$
we see from (\ref{fneu}) that $\zeta(u)$ has a pole of first order
in $\infty_k\in G_k$ for all $k=1,\dots,n-1$.\par Now let
$\zeta(u)\in G_k$ with some $k>1$. In this cases
$\zeta(u)\to\xi_k$ for $u\to\infty$. Again we obtain that $f$ has
a simple pole in $\infty_k\in G_k$ by (\ref{fneu}), because $q$
has only simple derivative zeros which are no zeros of $q$
(compare the remark at the beginning of section \ref{Riemann}).
\section{Blowing up and pulling back}
Let $r>0$ and $p_r(z)=r^np(z/r)$. If we start the considerations
of the preceding section with $p_r$ instead of $p$ we have to
replace the zeros $z_1,\dots,z_n$ of $p$ by $rz_1,\dots,rz_n$ and
the derivative zeros $\zeta(u)$ by $r\zeta(u)$ as well as $q(z)$
by $r^{n-1}q(z/r)$. The variation is then
$$Q_r(z,u):=r^nQ(z/r,u)=(z-ur)\cdot r^{n-1}q(z/r)=(z-u
r)\prod_{j=2}^n(z-z_jr).$$
\par Let $w_0$ be
an arbitrary complex number on the unit circle. If $r$ is large
enough me may provide that $|f(rw_0,r\zeta(w_0))|>1$, because of
the poles of $f$ in $\infty_k\in G_k$ for all $k=1,\dots,n-1$.
From (\ref{Qf}) we conclude that
$|r\zeta(w_0)-rw_0|>|r\zeta-rz_1|$ for all sufficiently large $r$
and all derivative zeros of $Q_r(.,w_0)$. This gives
$|\zeta(w_0)-w_0|>|\zeta-z_1|$. We define $Q(.,w_0)=:p^*$. If
$\zeta$ has been taken above as an essential derivative zero of
$p$ this says that $|\zeta(w_0)-w_0|>\rho(p,z_1)$ for all
derivative zeros of $p^*$. This gives $\rho(p^*,w_0)>\rho(p,z_1)$.
and theorem \ref{thm1} is proved.

\section{Proof of Sendov's conjecture}
The following Theorem has been proved by Goodman, Rahman, Ratti
\cite{grr} and independently by Schmeisser \cite{sch}.
\begin{thm} \label{grrs} Let $n\ge 2$ and $p\in\Pn$. If $p(1)=0$ then there is some
$\zeta$ with $p'(\zeta)=0$ and $|2\zeta-1|\le 1$. Moreover, there is such some
$\zeta$ in the open disk $|2z-1|< 1$ unless all derivative zeros of $p$
lie on the circle $|2z-1|= 1$ \end{thm}
The polynomial $z^n-1$ shows that
$\rho(p,z_1)\ge 1$ in order that $p$ is maximal with respect to its zero
$z_1$.
So the only point in the closed disk $|2z-1|\le 1$ where $p$ may have a
derivative zero is $0$. So we obtain from Theorem \ref{grrs} that $p'$ has a
single zero, located in the origin. Now  $p(1)=0$ implies $p(z)=z^n-1$.
We see that the only maximal polynomials in $\Pn$ are given by $z^n+a$ with
$|a|=1$. This
has been conjectured 1972 by Phelps and Rodriguez \cite{phro}.


\begin{minipage}[t]{14cm}
Universit\"at Oldenburg\\Fakult\"at V, Institut f\"ur Mathematik\\
Postfach 2503
\\ D-26111 Oldenburg\\ Bundesrepublik Deutschland\\ {\small
e-mail: schmieder@mathematik.uni-oldenburg.de}
\end {minipage}
\end{document}